\documentclass{amsart}
\usepackage{amssymb,amscd,amsthm,verbatim}
\usepackage{eucal}
\usepackage[dvips]{graphicx}
\usepackage{multirow}
\usepackage{fancyhdr}
\newtheorem{proposition}{Proposition}[section]

\newtheorem{prop}[proposition]{Proposition}
\newtheorem{thm}[proposition]{Theorem}
\newtheorem{cor}[proposition]{Corollary}
\newtheorem{lem}[proposition]{Lemma}

\newtheorem{defi}[proposition]{Definition}
\newtheorem{rem}[proposition]{Remark}

\begin{document}
\title{Method of Moments Estimation of Ornstein-Uhlenbeck Processes Driven by General L\'{e}vy Process}

\author{{\sc By Konstantinos Spiliopoulos}
\\{\sc Brown University, Providence, R.I. 02912, USA}\\ {\sc kspiliop@dam.brown.edu}}
\date{}

%

%

\maketitle

\begin{abstract}
Ornstein-Uhlenbeck processes driven by general L\'{e}vy process are
considered in this paper. We derive strongly consistent estimators
for the moments of the underlying L\'{e}vy process and for the mean
reverting parameter of a discretely observed L\'{e}vy driven
Ornstein-Uhlenbeck process. Moreover, we prove that the estimators
are asymptotically normal. We use ergodicity arguments. Finally, we
test the empirical performance of our estimators in a simulation
study and we fit the model to  real VIX data.
\end{abstract}

\thanks{AMS 2000 \textit{subject classifications.} Primary 62M05, 62G20; secondary 62P05.}

\section{Introduction}
Given a positive number $\lambda$ and a
time-homogeneous L\'{e}vy process $L$, the  Ornstein-Uhlenbeck (OU) process driven by $L$ is defined by
\begin{equation}
Y_{t}=e^{-\lambda t}Y_{0}+e^{-\lambda t}\int_{0}^{\lambda
t}e^{s}dL_{s}, \label{OU2}
\end{equation}
where $Y_{0}$ is assumed to be independent of $\{L_{t}\}_{t\geq
0}$. Following the terminology introduced by Barndorff-Nielsen and
Shephard in \cite{BNS1}, we shall call $L$ the background driving
L\'{e}vy process (BDLP). It is easy to see that (\ref{OU2}) is the
unique strong solution of the stochastic differential equation
\begin{eqnarray}
dY_{t}&=&-\lambda Y_{t}dt+dL_{\lambda t}.\label{OU1}
\end{eqnarray}
Under some regularity conditions on the L\'{e}vy measure of $L$ and
if $\lambda>0$, $Y$ admits a unique invariant distribution $F_{Y}$.
Owing to the scaling of the time index of $L$ in (\ref{OU1}) by
$\lambda$ (i.e. the term $L_{\lambda t}$), $F_{Y}$ is independent of
$\lambda$.

Let us suppose now that we have discrete-time observations
$Y_{0},Y_{h},\cdots$ $,Y_{(n-1)h}$ with $h>0$ from $\lbrace
Y_{t}\rbrace_{t\geq 0}$ as it is defined by (\ref{OU2}). The
objective here is to estimate the parameters of the model using
these discrete-time observations.  In particular, we are interested
in estimating $\lambda$ and moments of $L_{1}$. We derive strongly consistent method of moments
estimators and prove that they are asymptotically normal. In this paper, we
consider $\lambda$ and only the first two moments of $L_{1}$, i.e.
$\mathbb{E}L_{1}$ and $\mathbb{E}{L_{1}^{2}}$. However, the
methodology can be extended to higher moments as well  (see Remark \ref{Remark32}). Similar
methods of moments estimators have been used elsewhere as well (see
Valdivieso et al \cite{Valdivieso} for example), but without theoretical
justification. The proof of their consistency and asymptotic
normality is presented, according to the best knowledge of the
author, for the first time in the present work.

Our motivation for studying this problem comes from continuous
stochastic volatility models in financial mathematics.
Barndorff-Nielsen and Shephard (in \cite{BNS1}; see also
\cite{BNS2}) model stock price as a geometric Brownian motion and the
diffusion coefficient of this motion as an OU process that is driven
by a subordinator (a L\'{e}vy process that is nonnegative and
nondecreasing). Other continuous stochastic volatility models can be
found in Kl\"{u}ppelberg et al \cite{KLM1} and in Shephard \cite{Shephard}. Some papers that
consider statistical inference of these models are
Barndorff-Nielsen and Shephard \cite{BNS1},  Brockwell et al \cite{BrockwellDavisYang}, Haug et al \cite{HKLM}, Jongbloed et al \cite{JMV}.

The paper is organized as follows. Section 2 presents several known
results concerning OU processes and L\'{e}vy processes.  In section
3 we consider strongly consistent estimators of the first two
moments of $L_{1}$ and of $\lambda$ and we provide a methodology to
express any moment of the stationary distribution of $\{Y\}_{t\geq
0}$ in terms of the moments of $L_{1}$. In section 4, we prove that
these estimators are asymptotically normal. Section 5 discusses
modeling issues and simulation techniques and presents simulation
results for gamma OU process and inverse Gaussian OU process. In
section 6 we fit the model to real log(VIX) data and we argue that an OU model is a good candidate for modeling log(VIX). Finally, section 7 contains a summary and a discussion on future work.

We would like to mention here, that after completion of this work,
the author learned about the results in Jongbloed et al \cite{JMV}. In \cite{JMV},
the authors assume that $L$ is a subordinator. Let $F_{L}$ denote the L\'{e}vy measure of $L$, $Y$
the unique stationary solution to (\ref{OU2}) (which exists if $\int_{x>1}\log(x)F_{L}(dx)<\infty$ for example) and $F_{Y}$ its probability law.
The characteristic function of $Y$ is given by
\begin{displaymath}
\phi_{F_{Y}}(t):=\int e^{itx}F_{Y}(dx)=\exp(\int_{0}^{\infty}[e^{itx}-1]\frac{\kappa(x)}{x}dx),
\end{displaymath}
where $\kappa(x)=F_{L}(x,\infty)$. Hence, the stationary
distribution $F_{Y}$, of the OU process $Y$, is being determined
by the canonical function $\kappa(x)$. In \cite{JMV}, the authors
develop a nonparametric inference procedure for $\lambda$ and for
the canonical function $\kappa(x)$.  The results in the present
paper complement the results of \cite{JMV}.

\section{Assumptions and Preliminary Results}
Consider a probability space $(\Omega, \mathfrak{F}, P)$ equipped
with a filtration $\mathfrak{F}_{t}$.
\begin{defi}
 A one dimensional $\mathfrak{F}_{t}$ adapted
L\'{e}vy process is usually denoted by
$L_{t}=L_{t}(\omega)$, $t\geq 0$, $\omega\in \Omega$ and is a
stochastic process that satisfies the following:
\begin{enumerate}
\item{$L_{t}\in \mathfrak{F}_{t}$ for all $t\geq 0$.}
\item{$L_{0}=0$ a.s.}
\item{$L_{t}-L_{s}$ is independent of $\mathfrak{F}_{s}$ and has the same distribution as $L_{t-s}$.}
\item{It is a process continuous in probability.}
\end{enumerate}
\end{defi}
We assume that we are working with a c\`{a}dl\`{a}g L\'{e}vy
process (i.e. it is right continuous with left limits). It is well known that every L\'{e}vy process has such a modification.

Furthermore, if $F_{L}$ denotes the L\'{e}vy measure of $L_{1}$, we
will assume that there exist a constant $M>0$ such that
\begin{equation}
\int_{|x|>1}e^{v x}F_{L}(dx)<\infty, \hspace{0.5 cm} \textrm{ for
every } |v|\leq M. \label{MomentCondition}
\end{equation}
Condition (\ref{MomentCondition}) gaurantees that the moment
generating function $v\rightarrow \mathbb{E}e^{vL_{1}}$ exists at
least for $|v|\leq M$ (see Wolfe \cite{Wolfe} and Eberlein and
Raible \cite{EberleinRaible}).

We shall write
\begin{enumerate}
\item{$\mathbb{E}L_{1}=\mu$.}
\item{$\textrm{Var}(L_{1})=\sigma^{2}$.}
\end{enumerate}
Moreover, we shall assume that $Y_{0}$ is independent of $\{L_{t}\}_{t\geq 0}$ and that
\begin{equation}
Y_{0} \stackrel{\mathfrak{D}}{=} \int_{0}^{\infty}e^{-s}dL_{s}.\label{DistributionCondition}
\end{equation}
The integral on the right hand side of (\ref{DistributionCondition}) is well defined (see Sato \cite{Sato} for example).
The following proposition, which is a reformulation of Propositions 1 and 2 in Brockwell \cite{Brockwell}, characterizes the stationarity of the OU process $\{Y_{t}\}_{t\geq 0}$.
\begin{prop}
If $Y_{0}$ is independent of $\{L_{t}\}_{t\geq 0}$ and $EL_{1}^{2}<\infty$ then
$\{Y_{t}\}_{t\geq 0}$ is weakly stationary if and only if  $\lambda>0$ and $Y_{0}$ has the
same mean and variance as $\int_{0}^{\infty}e^{-s}dL_{s}$. If in addition
$Y_{0}$ has the same distribution as
$\int_{0}^{\infty}e^{-s}dL_{s}$, then $\{Y_{t}\}_{t\geq 0}$ is strictly
stationary and vice-versa.\label{Proposition22}
\end{prop}
In Masuda \cite{Masuda} now, the author proves, under mild regularity conditions, that the OU process $Y$ is strong Feller, its probability law has a smooth transition density, is ergodic and exponentially $\beta-$mixing (strong mixing). Before mentioning the results of \cite{Masuda} that we will use in the present paper, let us recall the definitions of a self-decomposable law on $\mathbb{R}$ and of $\beta-$mixing.
\begin{defi}
Let $\lambda$ be a positive number. Then, an infinitely divisible
distribution $F_{Y}$ is called $\lambda-$self-decomposable, if there
exists a random variable $X=X_{t,\lambda}$, such that, for each
$t\in\mathbb{R}_{+}$
\begin{displaymath}
\phi_{F_{Y}}(u)=\phi_{F_{Y}}(e^{-\lambda t}u)\phi_{F_{X}}(u), u\in \mathbb{R},
\end{displaymath}
where $\phi_{F_{Y}}(u)$ and $\phi_{F_{X}}(u)$ are the characteristic
functions corresponding to $F_{Y}$ and $F_{X}$ respectively. For the
sake of notational convenience we will just say that $F_{Y}$ is
called self-decomposable.
\end{defi}
If $\int_{|x|>1}\log(|x|)F_{L}(dx)<\infty$, then the class of all possible invariant distributions of $Y$ forms the class of all self-decomposable distributions $F_{Y}$ (see Sato \cite{Sato}). In particular, the latter is implied by (\ref{MomentCondition}).
\begin{defi}
For a stationary process $Y=\{Y_{t}\}_{t \geq 0}$ define the $\sigma$-algebras
$\mathfrak{F}_{1}=\mathfrak{F}_{(0,u)}=\sigma(\{Y_{v}\},0\leq v<u)$
and
$\mathfrak{F}_{2}=\mathfrak{F}_{[u+t,\infty)}=\sigma(\{Y_{v}\},v\geq
u+t)$. Then
\begin{enumerate}
\item{$Y$ is called $\beta-$mixing (or strong mixing) if:  $$\beta(t)=\sup_{A\in\mathfrak{F}_{1},B\in\mathfrak{F}_{2}}|P(A\cap B)-P(A)P(B)|\rightarrow 0 \textrm{ as } t\rightarrow \infty.$$}
\item{$Y$ is called $\beta-$mixing with exponential rate if for some $k>0$ and $a>0$: $$\beta(t)\leq k e^{-at} \textrm{ for } t\geq 0.$$}
\end{enumerate}
\end{defi}
The following theorem is Theorem 4.3 in Masuda \cite{Masuda} and discusses the mixing properties of $\{Y_{t}\}_{t\geq 0}$.
\begin{thm}
Let $\lambda>0$ and $\{Y_{t}\}_{t\geq 0}$ be the strictly stationary OU process given by (\ref{OU2})
 with self-decomposable marginal distribution $F_{Y}$. If we have that
$$\int_{\mathbb{R}}|x|^{p}F_{Y}dx<\infty$$
for some $p>0$, then there exists a constant $a>0$ such that $\beta(t)=O(e^{-at})$ as $t\rightarrow \infty$. In particular, $Y$ is ergodic.  \label{Theorem23}
\end{thm}

\section{Method of Moments Estimation}
We aim at estimation of the model parameters $\theta_{0}=(\mu,\sigma^{2},\lambda)$  from a sample of equally spaced observations from (\ref{OU2}) by matching moments and empirical autocorrelation function to their theoretical counterparts.

Proposition \ref{Proposition31} below relates the theoretical moments of $L_{1}$ with the theoretical
moments of  the stationary distribution $F_{Y}$ of $\{Y_{t}\}$.

\begin{prop}
Suppose that $\lbrace L_{t}\rbrace_{t\geq 0}$ is a L\'{e}vy process
such that $\mathbb{E}L_{1}=\mu<\infty$,
$\textrm{Var}{L_{1}}=\sigma^{2}<\infty$ and that
(\ref{MomentCondition}) holds. Let $M$ be the largest constant
satisfying (\ref{MomentCondition}) and assume that $\lambda<M$.
Then, the following are true
\begin{enumerate}
\item{$\mathbb{E}Y_{0}=\mu$}
\item{$\textrm{Var}{Y_{0}}=\frac{\sigma^{2}}{2}$}
\end{enumerate}\label{Proposition31}
\end{prop}

\begin{proof}
Let $\gamma(v)$ be the cumulant function of $L_{1}$, i.e.
\begin{equation}
\gamma(v)=\ln \mathbb{E}e^{vL_{1}}\label{CumulantOfL1}
\end{equation}
By the L\'{e}vy- Khinchine representation Theorem we get that $\gamma(v)$ has the form
\begin{equation}
\gamma(v)=b v+\frac{c}{2} v^{2}+\int_{\mathbb{R}}(e^{vx}-1-vx)F_{L}(dx), \label{CumulantOfL2}
\end{equation}
which is valid for $|v|\leq M$. Moreover, $\gamma$ is continuously
differentiable (see Lukacs \cite{Lukacs}).

Using the assumptions $\mathbb{E}L_{1}=\mu$ and $\textrm{Var}{L_{1}}=\sigma^{2}$ and relations (\ref{CumulantOfL1}) and (\ref{CumulantOfL2}), it is easy to see that $b=\mu$ and $c=\sigma^{2}-\int_{\mathbb{R}}x^{2}F_{L}(dx)$.

In order to calculate $\mathbb{E}Y_{0}$ and $\textrm{Var}{Y_{0}}$ we use the following formula:
\begin{equation}
\mathbb{E}e^{\int_{0}^{\infty}\lambda e^{-s}dL_{s}}=e^{\int_{0}^{\infty}\gamma(\lambda e^{-s})ds},\label{CumulantOfL3}
\end{equation}
which is valid since $\lambda<M$ (see Lemma 3.1 of  Eberlein and Raible \cite{EberleinRaible}).

Recall now that we have assumed $Y_{0}=\int_{0}^{\infty}e^{-s}dL_{s}$ in distribution. The latter and (\ref{CumulantOfL3}) imply that:
\begin{eqnarray}
\mathbb{E}Y_{0}&=&\frac{d}{d\lambda}\mathbb{E}e^{\int_{0}^{\infty}\lambda e^{-s}dL_{s}}|_{\lambda=0}=\nonumber\\
&=& \frac{d}{d\lambda}e^{\int_{0}^{\infty}\gamma(\lambda e^{-s})ds}|_{\lambda=0}=\nonumber\\
&=& \mu
\end{eqnarray}
In a similar way we get that $\mathbb{E}Y_{0}^{2}=\frac{\sigma^{2}}{2}+\mu^{2}$. This concludes the proof of the proposition.
\end{proof}
\begin{rem}
We would like to note here, that the proof of Proposition \ref{Proposition31} can be used for the calculation of higher moments of $Y_{0}$.\label{Remark32}
\end{rem}

It follows directly by (\ref{OU2})  that  the theoretical autocovariance and autocorrelation function of $Y_{t}$ are given by the formulas
\begin{enumerate}
\item{autocovariance: $\gamma(h)=cov(Y_{t+h},Y_{t})=\frac{\sigma^2}{2}e^{-\lambda h}$, for $h\in \mathbb{N}_{0}$.}
\item{autocorrelation: $\rho(h)=corr(Y_{t+h},Y_{t})=e^{-\lambda h}$, for $h\in \mathbb{N}_{0}$.}\label{TheoreticalAutoCorrelationOfY}
\end{enumerate}
On the other hand, the empirical moments, autocorrelation and
autocovariance function are given by the formulas below. Let $d\geq
0$ be fixed. Then, we have:
\begin{enumerate}
\item{Sample mean: $\bar{Y}_{\cdot}=\frac{1}{n}\sum_{i=1}^{n} Y_{i}$.}
\item{Sample variance: $\frac{1}{n}\sum_{i=1}^{n} (Y_{i}-\bar{Y}_{\cdot})^{2}$.}
\item{Sample autocovariance: $\hat{\gamma}_{n}=(\hat{\gamma}_{n}(0),\hat{\gamma}_{n}(1),\cdots,\hat{\gamma}_{n}(d))^{T}$ where for $h\in \lbrace 0,\cdots, d\rbrace$ we define $\hat{\gamma}_{n}(h)=\frac{1}{n}\sum_{i=1}^{n-h} (Y_{i+h}-\bar{Y}_{\cdot})(Y_{i}-\bar{Y}_{\cdot})$.}
\item{Sample autocorrelation: $\hat{\rho}_{n}=(\hat{\rho}_{n}(0),\hat{\rho}_{n}(1),\cdots,\hat{\rho}_{n}(d))^{T}$ where for $h\in \lbrace 0,\cdots, d\rbrace$ we define $\hat{\rho}_{n}(h)=\frac{\hat{\gamma}_{n}(h)}{\hat{\gamma}_{n}(0)}$.}\label{EmpiricalAutoCorrelationOfY}
\end{enumerate}
We have the following Theorem:
\begin{thm}
Let $\mu,\sigma^{2}, \gamma(\cdot), \hat{\gamma}_{n}(\cdot),\rho(\cdot)$ and $\hat{\rho}_{n}(\cdot)$ be defined as above. Then, the following statements are true
\begin{enumerate}
\item{$\bar{Y}_{\cdot}\stackrel{n\rightarrow\infty}{\longrightarrow} \mu$ almost surely}
\item{$\frac{1}{n}\sum_{i=1}^{n} (Y_{i}-\bar{Y}_{\cdot})^{2}\stackrel{n\rightarrow\infty}{\longrightarrow} \frac{\sigma^{2}}{2}$ almost surely}
\item{$(\hat{\gamma}_{n}(1),\cdots,\hat{\gamma}_{n}(d))\stackrel{n\rightarrow\infty}{\longrightarrow}(\gamma(1),\cdots,\gamma(d))$ almost surely}
\item{$(\hat{\rho}_{n}(1),\cdots,\hat{\rho}_{n}(d))\stackrel{n\rightarrow\infty}{\longrightarrow}(\rho(1),\cdots,\rho(d))$ almost surely}
\end{enumerate}\label{Theorem32}
\end{thm}
\begin{proof}
Due to our assumptions, the process $\{Y_{t}\}_{t\geq 0}$ is
strictly stationary. Moreover by Theorem \ref{Theorem23} it is also
$\beta-$mixing with exponential decaying rate. These two results
imply ergodicity of $\{Y_{t}\}_{t\geq 0}$. The latter together with
strict stationarity imply that empirical moments and sample
autocovariance functions are strongly consistent estimators of the
corresponding theoretical quantities Billingsley \cite{Billingsley}.
Then, the statement of the Theorem follows.
\end{proof}
For the mean reverting parameter $\lambda$ we have the following Lemma.
\begin{lem}
Let $K$ be a compact subset of $\mathbb{R}_{+}$
such that the true value of $\lambda$, say $\lambda_{o}$, belongs to
$K$ and let $\hat{\lambda}_{n}=\textrm{argmin}_{\lambda \in
K}\sum_{h=1}^{d}(\hat{\rho}_{n}(h)-e^{-\lambda h})^{2}$. Then $\hat{\lambda}_{n}$ exists, is locally unique and
\begin{equation}
\hat{\lambda}_{n}\stackrel{n\rightarrow\infty}{\longrightarrow} \lambda_{o} \textrm{
almost surely}.
\end{equation}\label{Lemma33}
\end{lem}
\begin{proof}
Consider the functions $\Delta_{n}(\lambda)=\sum_{h=1}^{d}(\hat{\rho}_{n}(h)-\rho_{\lambda}(h))^{2}$ and $\Delta_{0}(\lambda)=\sum_{h=1}^{d}(\rho_{\lambda_{o}}(h)-\rho_{\lambda}(h))^{2}$
where $\rho_{\lambda}(h)=e^{-\lambda h}$. Theorem \ref{Theorem32} implies that for all $\lambda\in K$:
\begin{displaymath}
\Delta_{n}(\lambda)\stackrel{n\rightarrow\infty}{\longrightarrow}\Delta_{0}(\lambda) \hspace{0.2cm} \textrm{almost surely}.
\end{displaymath}
By Theorem II.1 in Andresen and Grill \cite{AndresenGrill} we have
\begin{displaymath}
\sup_{\lambda\in K}|\Delta_{n}(\lambda)-\Delta_{0}(\lambda)|\stackrel{n\rightarrow\infty}{\longrightarrow}0 \hspace{0.2cm} \textrm{almost surely}.
\end{displaymath}
Observe now that $\Delta_{0}(\lambda)$ is a sum of nonnegative
terms. It becomes zero if and only if $\lambda=\lambda_{0}$. Hence,
$\Delta_{0}(\lambda)$ has a unique minimum at $\lambda=\lambda_{0}$
which is equal to zero. We get
\begin{displaymath}
\Delta_{n}(\lambda_{0})\stackrel{n\rightarrow\infty}{\longrightarrow}0 \hspace{0.2cm}  \textrm{almost surely}.
\end{displaymath}
Furthermore, for $n$ finite we have that $0\leq \Delta_{n}(\hat{\lambda}_{n})\leq \Delta_{n}(\lambda_{0})$. Therefore we get
\begin{displaymath}
\Delta_{n}(\hat{\lambda}_{n})\stackrel{n\rightarrow\infty}{\longrightarrow}0 \hspace{0.2cm} \textrm{almost surely}.
\end{displaymath}
Moreover, we have
\begin{eqnarray}
|\Delta_{n}(\hat{\lambda}_{n})-\Delta_{0}(\hat{\lambda}_{n})|&=&
|\sum_{h=1}^{d}[\hat{\rho}_{n}^{2}(h)-\rho_{\lambda_{0}}^{2}(h)+2\rho_{\hat{\lambda}_{n}}(h)
(\rho_{\lambda_{0}}(h)-\hat{\rho}_{n}(h))]|\leq\nonumber\\
&\leq&\sum_{h=1}^{d}[|\hat{\rho}_{n}(h)|+|\rho_{\lambda_{0}}(h)|+2|\rho_{\hat{\lambda}_{n}}(h)|]|\rho_{\lambda_{0}}(h)-\hat{\rho}_{n}(h)|\leq\nonumber\\
&\leq& 4\sum_{h=1}^{d}|\rho_{\lambda_{0}}(h)-\hat{\rho}_{n}(h)|\stackrel{n\rightarrow\infty}{\longrightarrow} 0 \hspace{0.2cm} \textrm{almost surely.}\nonumber
\end{eqnarray}
Here, we used the relation $|\hat{\rho}_{n}(h)|\leq 1$ which follows
immediately from Cauchy-Schwarz inequality. The above imply that
\begin{displaymath}
\Delta_{0}(\hat{\lambda}_{n})\stackrel{n\rightarrow\infty}{\longrightarrow}0 \hspace{0.2cm} \textrm{almost surely}.
\end{displaymath}
But, $\lambda_{0}$ is the unique minimum of $\Delta_{0}(\lambda)$ and it satisfies $\Delta_{0}(\lambda_{0})=0$. Thus,
\begin{displaymath}
\Delta_{0}(\hat{\lambda}_{n})\stackrel{n\rightarrow\infty}{\longrightarrow}\Delta_{0}(\lambda_{0})=0 \hspace{0.2cm} \textrm{almost surely}.
\end{displaymath}
Hence, we easily conclude (Corollary II.2 in \cite{AndresenGrill}) that $\hat{\lambda}_{n}$ is locally uniquely determined and that
\begin{displaymath}
\hat{\lambda}_{n}\stackrel{n\rightarrow\infty}{\longrightarrow} \lambda_{o} \textrm{
almost surely}.
\end{displaymath}
\end{proof}
\begin{rem}
Theorem \ref{Theorem32} and Lemma \ref{Lemma33} give us two strongly
consistent estimators for $\lambda$. The first one is
$\hat{\lambda}_{1,n}=-\log(\hat{\rho}_{n}(1))$ and the second one is
$\hat{\lambda}_{2,n}=\textrm{argmin}_{\lambda}\sum_{h=1}^{d}(\hat{\rho}_{n}(h)-e^{-\lambda
h})^{2}$. One could use, for example, $\hat{\lambda}_{1,n}$  as an
initial value to an algorithm that calculates $\hat{\lambda}_{2,n}$.
\end{rem}
Summarizing, we have that $\hat{\mu}_{n}$, $\hat{\sigma}^{2}_{n}$
and $\hat{\lambda}_{1,n}, \hat{\lambda}_{2,n}$ are strongly
consistent estimators of $\mu,\sigma^{2}$ and $\lambda$
respectively, where:
\begin{eqnarray}
\hat{\mu}_{n}&=&\frac{1}{n}\sum_{i=1}^{n} Y_{i}\nonumber\\
\hat{\sigma}^{2}_{n}&=&2\frac{1}{n}\sum_{i=1}^{n} (Y_{i}-\hat{\mu}_{n})^{2} \label{StrongConsistentEstimators}\\
\hat{\lambda}_{1,n}&=&-\log(\hat{\rho}_{n}(1))\nonumber\\
\hat{\lambda}_{2,n}&=&\textrm{argmin}_{\lambda}\sum_{h=1}^{d}(\hat{\rho}_{n}(h)-e^{-\lambda h})^{2}.\nonumber
\end{eqnarray}
\begin{rem}
For a stationary model, the parameter $\lambda$ has to be positive.
However, if we compute $\hat{\lambda}_{2,n}$ as the unrestricted
minimum $\hat{\lambda}_{2,n}=\textrm{argmin}_{\lambda \in
\mathbb{R}_{+}}\sum_{h=1}^{d}(\hat{\rho}_{n}(h)-e^{-\lambda h})^{2}$
we may end up with a negative estimator $\hat{\lambda}_{n}$. In this
case, we define the estimator of $\lambda$ to be zero and we take
this as an indication that the data is not stationary.
\end{rem}

\section{Asymptotic Properties of the Moment Estimators}

In this section we prove that the estimators defined by (\ref{StrongConsistentEstimators}) are asymptotically normal.

If $\beta$ is a vector, then we define by $\beta^{T}$ its transpose.
We begin with the following central limit theorem.

\begin{thm}
Let us assume that there exists a $\delta>0$ such that $\mathbb{E}Y_{0}^{4+\delta}<\infty$.
Define
\begin{eqnarray}
\hat{\psi}_{n}&=&(\hat{\mu}_{n},\hat{\gamma}_{n}(0),\hat{\gamma}_{n}(1),\cdots,\hat{\gamma}_{n}(d))^{T}\nonumber\\
\psi_{o}&=&(\mu,\gamma(0),\gamma(1),\cdots,\gamma(d))^{T}\nonumber\\
\Sigma&=&[\sigma_{k,l}]_{k,l=1}^{d+2}\textrm{ with elements
}\nonumber\\
&
&\sigma_{k,l}=cov(Z_{1}^{k},Z_{1}^{l})+2\sum_{i=1}^{\infty}cov(Z_{1}^{k},Z_{i+1}^{l})\textrm{ where }\nonumber\\
Z_{i}&=&(Y_{i},(Y_{i}-\mu)^{2},
(Y_{i+1}-\mu)(Y_{i}-\mu),\cdots,(Y_{i+d}-\mu)(Y_{i}-\mu))^{T}\nonumber
\end{eqnarray}
Then, the following holds:
\begin{equation}
\sqrt
n(\hat{\psi}_{n}-\psi_{o})\stackrel{\mathfrak{D}}{\longrightarrow} N(0, \Sigma)\label{CLT1}
\end{equation}
where $N(0,\Sigma)$ is the multivariate normal distribution with mean $0$ and variance-covariance matrix $\Sigma$.\label{Theorem41}
\end{thm}
\begin{proof} The proof of this theorem is similar to the proof of Proposition 3.7 of Haug et al \cite{HKLM}. Let us define
\begin{enumerate}
\item{$\gamma_{n}^{*}(h)=\frac{1}{n}\sum_{i=1}^{n} (Y_{i+h}-\mu)(Y_{i}-\mu), h\in\lbrace 0,\cdots, d\rbrace.$}
\item{$\gamma_{n}^{*}=(\gamma_{n}^{*}(0),\cdots,\gamma_{n}^{*}(d))^{T}$.}
\end{enumerate}
We first prove that (\ref{CLT1}) is true with $\hat{\psi}_{n}^{*}=(\hat{\mu}_{n},\hat{\gamma}_{n}^{*}(0),\hat{\gamma}_{n}^{*}(1),\cdots,\hat{\gamma}_{n}^{*}(d))^{T}$ in place of $\hat{\psi}_{n}$.

By the well known Cramer-Wold device, it is sufficient to prove that for every $\beta\in \mathbb{R}^{d+2}$ such that $\beta^{T}\Sigma\beta>0$ we have
\begin{equation}
\sqrt{n}(\frac{1}{n}\sum_{i=1}^{n}\beta^{T}Z_{i}-\beta^{T}\psi_{0})\stackrel{\mathfrak{D}}{\longrightarrow}N(0, \beta^{T}\Sigma\beta).\label{CLT2}
\end{equation}
It is well known (see \cite{Billingsley} for example) that strong mixing and the corresponding decaying rate are preserved under linear transformations. Thus, the sequence $\lbrace\beta^{T}Z_{i}\rbrace$ is strong mixing with exponential decaying rate. Since, by assumption $\mathbb{E}|Z_{1}|^{2+\epsilon}$ for some $\epsilon>0$, the central limit theorem for strong mixing processes is applicable (Theorem 7.3.1 in Ethier and Kurtz \cite{EK}). Hence, we have as $n\rightarrow \infty$ that
\begin{equation}
\sqrt{n}(\frac{1}{n}\sum_{i=1}^{n}\beta^{T}Z_{i}-\beta^{T}\psi_{0})\stackrel{\mathfrak{D}}{\longrightarrow} N(0, \tilde{\sigma}^{2}).\label{CLT3}
\end{equation}
But, we easily see that $\tilde{\sigma}^{2}=var(\beta^{T}Z_{1})+2\sum_{i=1}^{\infty}cov(\beta^{T}Z_{1},\beta^{T}Z_{i+1})=\beta^{T}\Sigma\beta$. So (\ref{CLT2}) holds.

Now recall that by Theorem \ref{Theorem32} we have
\begin{equation}
\hat{\psi}_{n}\stackrel{n\rightarrow\infty}{\longrightarrow}\psi_{0}\hspace{0.2cm} \textrm{ almost surely}.
\end{equation}
Following the proof of proposition 7.3.4 of Brockwell and Davis \cite{BrockwellDavis} we get
\begin{equation}
\sqrt{n}(\frac{1}{n}\sum_{i=1}^{n}\beta^{T}Z_{i}-\beta^{T}\hat{\psi}_{n})\stackrel{n\rightarrow\infty}{\longrightarrow} 0 \hspace{0.2cm} \textrm{in probability}.\label{CLT4}
\end{equation}
Therefore, $\hat{\psi}_{n}$ has the same asymptotic behavior as $\hat{\psi}_{n}^{*}$. The latter and  (\ref{CLT2}) imply the Theorem.
\end{proof}
\begin{cor}
Let the conditions of Theorem \ref{Theorem41} hold. Then we have
\begin{displaymath}
\sqrt{n}(\hat{\rho}_{n}-\rho)\stackrel{\mathfrak{D}}{\longrightarrow} N(0, \Sigma_{\rho}).\label{CLT5}
\end{displaymath}
\end{cor}
\begin{proof}
It follows directly by Theorem \ref{Theorem41} and delta method (Theorem 3.1 in A.W.van der Vaart \cite{VanDerVaart}).
\end{proof}
Finally, we prove central limit theorem for $\hat{\theta}_{n}=(\hat{\mu}_{n},\hat{\sigma}^{2}_{n},\hat{\lambda}_{2,n})^{T}$. Let us denote $\sigma^{2}_{Y}=\gamma(0)=\textrm{Var}(Y_{0})=\frac{\sigma^{2}}{2}$ and define the following mappings.

\begin{equation}
G:\mathbb{R}\times [0,\infty)^{2}\longrightarrow \mathbb{R}\times
[0,\infty)^{2}: G(\mu,\sigma^{2}_{Y},\lambda)=
\begin{cases}
(\mu,2\sigma^{2}_{Y},\lambda), & \lambda>0 \cr
   (\mu,2\sigma^{2}_{Y},0), & \lambda\leq 0. \cr
   \end{cases}
\end{equation}
\begin{equation}
F:\mathbb{R}^{d+1}_{+}\longrightarrow \mathbb{R}_{+}:
F(\hat{\rho})=\textrm{argmin}_{\lambda}\sum_{h=0}^{d}(\hat{\rho}_{n}(h)-e^{-\lambda h})^{2}=\hat{\lambda}_{2,n}
\end{equation}
and $H$ as follows:
\begin{equation}
H:\mathbb{R}^{d+2}\longrightarrow \mathbb{R}\times [0,\infty)^{2}:
H(\mu,\gamma^{T})=G(\mu,\sigma^{2}_{Y},F(\rho)),
\end{equation}
where $\rho(h)=\frac{\gamma(h)}{\gamma(0)}$ for $h=0,\cdots,d$.
\begin{thm}
Let the conditions of Theorem \ref{Theorem41} be satisfied. Let us define:
\begin{eqnarray}
\hat{\theta}_{n}&=&(\hat{\mu}_{n},\hat{\sigma}^{2}_{n},\hat{\lambda}_{2,n})^{T}\nonumber\\
\theta_{o}&=&(\mu,\sigma^{2},\lambda)^{T}.\nonumber
\end{eqnarray}
Then the following holds:
\begin{equation}
\sqrt
n(\hat{\theta}_{n}-\theta_{o})\stackrel{\mathfrak{D}}{\longrightarrow}[\frac{\partial
H(\mu,\gamma^{T})}{\partial
(\mu,\gamma^{T})}] N(0, \Sigma)
\end{equation}
\end{thm}
\begin{proof}
It follows directly by Theorem \ref{Theorem41} and delta method applied to the differentiable map $H$.
\end{proof}

\section{Modeling and Simulation}
In this section we discuss modeling issues of L\'{e}vy driven OU
processes and present some simulation results for a gamma OU process
and an inverse Gaussian OU process. We use the simulated
data to test the performance of our estimators. However, first we mention the necessary ingredients for the simulation process.

A very important ingredient in modeling of L\'{e}vy driven OU processes is the connection between the L\'{e}vy density of the stationary distribution of $Y$ to the L\'{e}vy density of the probability law of $L_{1}$. In particular we have the following proposition.

\begin{prop} Assume that the L\'{e}vy density of $Y$, $\nu_{Y}(x)$, is differentiable and denote the L\'{e}vy density of the probability law of $L_{1}$ by $\nu_{L}(x)$. Then the following relation holds.
\begin{equation}
\nu_{L}(x)=-\nu_{Y}(x)-x\nu'_{Y}(x).\label{RelationForLevyDensities}
\end{equation}\label{Proposition51}
\end{prop}
\begin{proof}
It follows directly by the fact that the stationary solution, $Y$, to $(\ref{OU1})$ satisfies
\begin{displaymath}
Y\stackrel{\mathfrak{D}}{=}\int_{0}^{\infty}e^{-\lambda s}dL(\lambda
s).
\end{displaymath}
See \cite{BNS1} and \cite{BNS2} for more details.
\end{proof}
Hence, given $\nu_{L}(x)$ we can find $\nu_{Y}(x)$ and vice-versa. One can specify the law of the one dimensional marginal distribution of the OU process $Y$ and work out the density of the BDLP, $L_{1}$. One can also go the other way and model through the BDLP. Of course, there are constraints on valid BDLP's which must be satisfied. In particular, if
\begin{displaymath}
\int_{\mathbb{R}}\min\{1,x^{2}\}\nu_{L}(x)dx
\end{displaymath}
then $\nu_{L}(x)$ is the density of a L\'{e}vy jump process $L$ and there exists an OU process $Y$ such that $L$ is the BDLP of $Y$.
A very good survey on the relation between several distributions of $Y$ and $L$ is Barndorff-Nielsen and Shephard \cite{BNS2} (see also Schoutens \cite{Schoutens}).

Another important ingredient in simulations is the infinite series representation of L\'{e}vy integrals (see Rosinski \cite{Rosinski}). For simplicity, we restrict attention to L\'{e}vy processes, $L$, that are subordinators, i.e. they are nonnegative and nondecreasing. It is easy to see that subordinators have no Gaussian component, nonnegative drift and a L\'{e}vy measure that is zero on the negative half-line. If $Y$ models stochastic volatility then it has to be positive and such a choice of the BDLP guarantees that.

Let us denote by $\Gamma_{L}^{+}$ the tail mass function of $\nu_{L}$, i.e.
\begin{equation}
\Gamma_{L}^{+}(x)=\int_{x}^{\infty}\nu_{L}(y)dy\label{ITMF1}
\end{equation}
and by $\Gamma_{L}^{-1}$ the generalized inverse function of $\Gamma_{L}^{+}$, i.e.
\begin{equation}
\Gamma_{L}^{-1}(x)=\inf\lbrace y>0:\Gamma_{L}^{+}(y)\leq x\rbrace \label{ITMF2}.
\end{equation}
In order to simulate from (\ref{OU2}) we need to be able to simulate from $e^{-\lambda t}\int_{0}^{\lambda
t}e^{s}dL_{s}$. The key result here is the following infinite series representation of this type of integrals (Rosinski \cite{Rosinski}):
\begin{prop}
Consider a subordinator $L$ with positive increments. Let $f$ be a positive and integrable function on $[0,T]$. Then
\begin{equation}
\int_{0}^{T}f(s)dL_{s}=\sum_{i=1}^{\infty}\Gamma_{L}^{-1}(\alpha_{i}/T)f(T r_{i}),\label{SeriesRepresentation}
\end{equation}
where the equality is understood in distributional sense, $\lbrace \alpha_{i}\rbrace$ and $\lbrace r_{i}\rbrace$ are two independent sequences of random variables such that $r_{i}$ are independent copies of a uniform random variable in $[0,1]$ and $\lbrace\alpha_{i}\rbrace$ is a strictly increasing sequence of arrival times of a Poisson process with intensity $1$.
\end{prop}
\begin{rem}
We note here that the convergence of the series (\ref{SeriesRepresentation}) is often quite slow.\label{Remark53}
\end{rem}
Using (\ref{SeriesRepresentation}) we can then simulate a L\'{e}vy driven OU process. In particular, if $\Delta$ denotes the time step, we will use the identity
\begin{eqnarray}
Y_{t+\Delta}&=&e^{-\lambda \Delta}(Y_{t}+e^{-\lambda t}\int_{ t}^{
t+\Delta}e^{\lambda s}dL_{\lambda s})\nonumber\\
&=&e^{-\lambda \Delta}(Y_{t}+\int_{0}^{
\Delta}e^{\lambda s}dL_{\lambda s}). \label{OUsimulation}
\end{eqnarray}
Let us demonstrate the validity of our estimators modeling through
the BDLP. We consider two cases: $(a)$ when $Y_{0}\sim
\textrm{Gamma}(a,b)$ and $(b)$ when $Y_{0}\sim \textrm{IG}(a,b)$,
where IG stands for inverse Gaussian.

Regarding the $\lambda$ parameter, we recall our estimators:
$\hat{\lambda}_{1,n}=-\frac{\log(\hat{\rho}_{n}(1))}{\Delta}$ and
$\hat{\lambda}_{2,n}=\textrm{argmin}_{\lambda}\sum_{h=1}^{d}(\hat{\rho}_{n}(h)-e^{-\lambda
h \Delta})^{2}$. In (\ref{StrongConsistentEstimators}) we defined
$\hat{\lambda}_{1,n}$ and $\hat{\lambda}_{2,n}$ for $\Delta=1$, but
of course one can generalize them to any $\Delta>0$.

\subsection{Gamma OU model}
Assume that the driving L\'{e}vy process $L$ is a compound Poisson process and in particular, that $L_{t}=\sum_{n=1}^{N_{t}}x_{n}$ where $N_{t}$ is Poisson with intensity parameter $a$ and $x_{n}$ are independent identically distributed $\textrm{Gamma}(1,b)$ random variables. Using (\ref{RelationForLevyDensities}) we  get that $Y_{0}\sim \textrm{Gamma}(a,b)$.  It is known (see \cite{BNS1}) that in this case
\begin{displaymath}
\Gamma_{L}^{-1}(x)=\max\lbrace 0,-\frac{1}{b}\log(\frac{x}{a})\rbrace.
\end{displaymath}
Using this and equations (\ref{SeriesRepresentation}) and (\ref{OUsimulation}) we can easily simulate from a $\textrm{Gamma}(a,b)$-OU process. We also need to know how the parameters $\mu$ and $\sigma^2$ relate to $a$ and $b$. Since $\mathbb{E}L_{1}=\mu$ and $\textrm{Var}(L_{1})=\sigma^{2}$ implies $\mathbb{E}Y_{0}=\mu$ and $\textrm{Var}(Y_{0})=\frac{\sigma^{2}}{2}$, we have that $a=2\frac{\mu^{2}}{\sigma^{2}}$ and $b=2\frac{\mu}{\sigma^{2}}$.

\vspace{0.1cm}

We simulated 100 independent paths of a gamma OU process of 1000
observations each, with time step $\Delta=0.1$, using
(\ref{OUsimulation}). We chose $\mu=2$, $\sigma^{2}=0.25$ and in
order to capture possible different behaviors of the intensity
parameter we chose two different values for $\lambda,$ $0.5$ and
$5$.

\vspace{0.1cm}

Tables I and II, summarize the results for $\theta_{0}=(2,0.25,0.5)$
and for $\theta_{0}=(2,0.25,5)$ respectively.

\begin{table}[!h]
\begin{center}
\begin{tabular}{|c|c|c|c|c|}
\hline
   True Values & Est. Values  & Sample Std. Error & Comments   \\
  \hline $\mu=2$ & $1.995458$ & $0.0702198$ &  - \\
  \hline $\sigma^{2}=0.25$ & $0.2350207$ & $ 0.05352894$ & - \\
  \hline $\lambda=0.5$ & $0.566116$ & $0.1126439$ & $\hat{\lambda}_{n}=\hat{\lambda}_{1,n}$ \\
  \hline $\lambda=0.5$ & $0.5879571$ & $ 0.1441501$ & $\hat{\lambda}_{n}=\hat{\lambda}_{2,n}$\\
  \hline
\end{tabular}
\end{center}
\caption{$\theta_{0}=(2,0.25,0.5)$}
\end{table}
\begin{table}[!h]
\begin{center}
\begin{tabular}{|c|c|c|c|c|}
\hline
   True Values & Est. Values  & Sample Std. Error & Comments   \\
  \hline $\mu=2$ & $2.003799$ & $0.02094129$ &  - \\
  \hline $\sigma^{2}=0.25$ & $ 0.2473567$ & $ 0.01608991$ & - \\
  \hline $\lambda=5$ & $5.12962$ & $0.4463517$ & $\hat{\lambda}_{n}=\hat{\lambda}_{1,n}$ \\
  \hline $\lambda=5$ & $5.186585$ & $0.5898125$ & $\hat{\lambda}_{n}=\hat{\lambda}_{2,n}$ \\
  \hline
\end{tabular}
\end{center}
\caption{$\theta_{0}=(2,0.25,5)$}
\end{table}

\newpage

\subsection{Inverse Gaussian OU model}
It is well known that if $Y_{0}\sim \textrm{IG}(a,b)$, then the
L\'{e}vy density of $Y$ is
\begin{equation}
\nu_{Y}(x)=\frac{1}{\sqrt{2\pi}}ax^{-3/2}e^{-\frac{1}{2}b^{2}x}.
\end{equation}
Consider now the Lambert-W function, $L_{w}(\cdot)$, which satisfies $L_{w}(x)e^{L_{w}(x)} =x$. As it is also shown in Gander and Stephens \cite{GanderStephens}, equations (\ref{RelationForLevyDensities}), (\ref{ITMF1}) and (\ref{ITMF2}) imply that the inverse tail mass function of the BDLP of an IG(a,b)-OU process is given by
\begin{equation}
\Gamma_{L}^{-1}(x)=\frac{1}{b^2}L_{w}(\frac{a^{2}b^{2}}{2\pi x^{2}})\label{ITMF3}.
\end{equation}
Using the latter and equations (\ref{SeriesRepresentation}) and (\ref{OUsimulation}) we can easily simulate an IG(a,b)-OU process. We also need to know how the parameters $\mu$ and $\sigma^2$ relate to $a$ and $b$. Since $EL_{1}=\mu$ and $Var(L_{1})=\sigma^{2}$ implies $EY_{0}=\mu$ and $Var(Y_{0})=\frac{\sigma^{2}}{2}$, we have that $a=\mu\sqrt{\frac{2\mu}{\sigma^{2}}}$ and $b=\sqrt{\frac{2\mu}{\sigma^{2}}}$.

\vspace{0.1cm}

We simulated 100 independent paths of an IG-OU process of 1000
observations each with time step $\Delta=0.1$, using
(\ref{OUsimulation}). As before, we chose $\mu=2$, $\sigma^{2}=0.25$
and two different values for $\lambda,$ $0.5$ and $5$.

Tables III and IV, summarize the results for $\theta_{0}=(2,0.25,0.5)$
and for $\theta_{0}=(2,0.25,5)$ respectively.

\begin{table}[!h]
\begin{center}
\begin{tabular}{|c|c|c|c|c|}
\hline
   True Values & Est. Values  & Sample Std. Error & Comments   \\
  \hline $\mu=2$ & $ 1.986862 $ & $ 0.06476202 $ &  - \\
  \hline $\sigma^{2}=0.25$ & $ 0.2331244 $ & $ 0.05235387 $ & - \\
  \hline $\lambda=0.5$ & $ 0.5581237$ & $ 0.1128397$ & $\hat{\lambda}_{n}=\hat{\lambda}_{1,n}$ \\
  \hline $\lambda=0.5$ & $ 0.6050457$ & $ 0.1376689$ & $\hat{\lambda}_{n}=\hat{\lambda}_{2,n}$\\
  \hline
\end{tabular}
\end{center}
\caption{$\theta_{0}=(2,0.25,0.5)$}
\end{table}

\vspace{0.1cm}

\begin{table}[!h]
\begin{center}
\caption{$\theta_{0}=(2,0.25,5)$}
\begin{tabular}{|c|c|c|c|c|}
\hline
   True Values & Est. Values  & Sample Std. Error & Comments   \\
  \hline $\mu=2$ & $ 1.955288$ & $ 0.03107831$ &  - \\
  \hline $\sigma^{2}=0.25$ & $  0.2452349$ & $ 0.01750871$ & - \\
  \hline $\lambda=5$ & $5.05211$ & $ 0.4262788$ & $\hat{\lambda}_{n}=\hat{\lambda}_{1,n}$ \\
  \hline $\lambda=5$ & $ 5.158421$ & $ 0.659508$ & $\hat{\lambda}_{n}=\hat{\lambda}_{2,n}$ \\
  \hline
\end{tabular}
\end{center}
\end{table}

\newpage


\section{Real Data Analysis}
In 1993, the Chicago Board Options Exchange (CBOE) introduced the
CBOE volatility index, VIX, and it quickly became a popular measure
for stock market volatility. In 2003, the VIX methodology was
updated (see www.cboe.com for more details on the old and new VIX
methodology).   VIX measures the implied volatility of S\&P $500$
index options and it provides a minute-by-minute  snapshot of the
markets expectancy of volatility over the next $30$ calendar days.

We fitted the gamma OU model and the IG OU model to daily log
opening values of the VIX for the year 2004 (VIX values are
calculated using the new methodology). The data are taken from
www.cboe.com. 
We use the values from $1/2/2004$ till $9/30/2004$ for the
calibration of the model (in total $189$ data points) and the values
from $10/1/2004$ till $11/30/2004$ (in total $41$ data points) for
testing the model.



Table V summarizes the estimators, given by (\ref{StrongConsistentEstimators}), of the parameters of the model. We used $\hat{\lambda}_{2,n}$ to estimate $\lambda$.
\begin{table}[!h]
\begin{center}
\begin{tabular}{|c|c|c|c|}
\hline
    parameter & $\mu$& $\sigma^2$ & $\lambda$      \\
   \hline estimated value & $2.781769$ & $0.01919740$ & $  0.1767250$\\
  \hline
\end{tabular}
\end{center}
\caption{Estimated values for the parameters.}
\end{table}

\vspace{0.2cm}

In Figure 1, we see the first 10 lags of the
empirical autocorrelation function of the $\log(\textrm{VIX})$ for
$1/2/2004$ till $9/30/2004$ versus the theoretical autocorrelation
function of the OU model with $\lambda=0.1767250$, i.e.
$\rho(h)=e^{-0.1767250 h}$.

\begin{figure}[!h]
\begin{center}
\includegraphics[width=11 cm, height=4.5 cm]{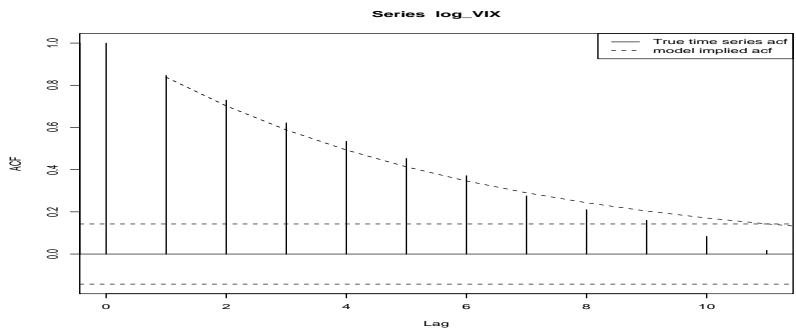}
\caption{ True time series acf versus the model implied acf with
$\lambda=0.1767250$.}
\end{center}
\end{figure}

As we saw before, the autocorrelation function of an OU model is
exponentially decreasing, i.e. it has the form  $e^{-\lambda h}$.
Figure 1 shows that  $e^{-0.1767250 h}$  approximates sufficiently
well the empirical autocorrelation function of $\log(\textrm{VIX})$
for $1/2/2004$ till $9/30/2004$, which is also exponentially
decreasing. Hence, we conclude that an OU model is a good candidate
for describing this data set.

To investigate the model fit, we performed a Ljung-Box test for the
squared residuals. We used the estimated values from Table V and since our data is daily
opening values we chose $\Delta=1$ for the time step. The test statistic used 10 lags of the empirical
autocorrelation function. The gamma OU model performed better than
the the IG-OU model. The null hypothesis was not rejected at the
0.05 level and the $p-$value was quite high, $0.6165$. In Figure 2
we see the empirical autocorrelation function of the residuals of $\log(\textrm{VIX})$ and in Figure 3 we see the actual residuals of the gamma OU model.

In Figure 4 we see in one figure: the actual time series from
$10/1/2004$ till $11/30/2004$, the one step ahead predicted time
series and $95\%$ bootstrap upper and lower confidence bounds of the
one step ahead predicted time series. In order to create the one
step ahead predicted time series we averaged over $50$ paths. We observe that the real time series (solid line) is most of the
time within the 95\% bootstrap upper and lower confidence bounds of
the one step ahead predicted time series (dotted lines), with very
few exceptions.

\begin{figure}[!h]
\begin{center}
\includegraphics[width=11 cm, height=4.0 cm]{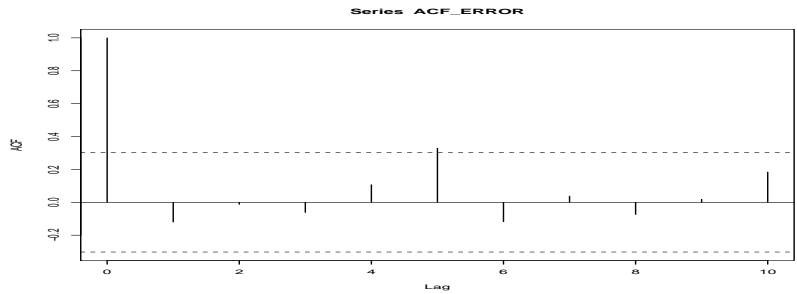}
\caption{ Empirical acf for the residuals of $\log(\textrm{VIX})$.}
\end{center}
\end{figure}

\begin{figure}[!h]
\begin{center}
\includegraphics[width=11 cm, height=4.0 cm]{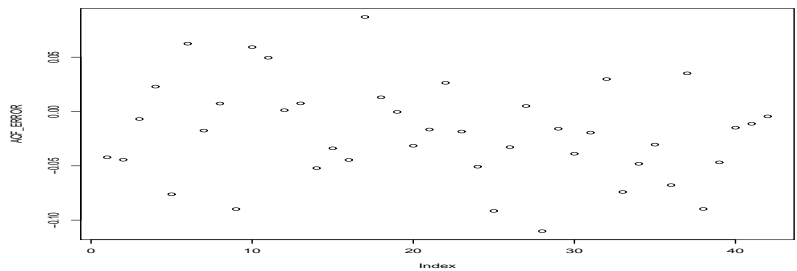}
\caption{The actual residuals.}
\end{center}
\end{figure}

\begin{figure}[!h]
\begin{center}
\includegraphics[width=11 cm, height=4 cm]{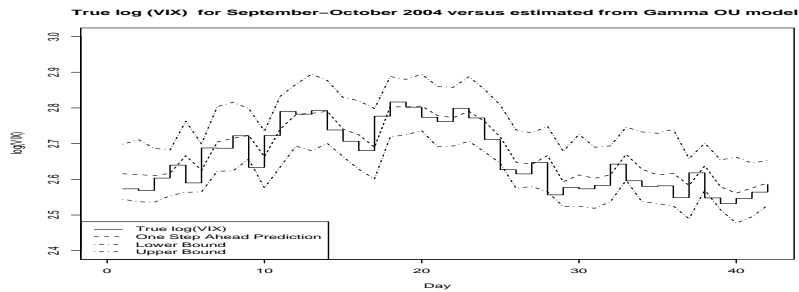}
\caption{Actual time series versus one step ahead predicted time series.}
\end{center}
\end{figure}

\newpage

\section{Discussion And Future Work}
In this paper, we consider an Ornstein-Uhlenbeck process driven by a
general L\'{e}vy process. We derive strong consistent estimators for
the parameters of the model and we prove that they are
asymptotically normal. Using simulated data, we show that the
estimators perform well at least for a gamma OU model and for an IG-OU model. Lastly we fit
the model to real data and we see that a L\'{e}vy driven OU model is a good
candidate for describing $\log(\textrm{VIX})$.



There are some interesting extensions to the model
studied in this paper. One such extension is a coupled two
dimensional OU process driven by a two dimensional L\'{e}vy process.
This model is important for financial applications, since it could
be used to model log of the price and stochastic volatility
simultaneously. An interesting question for financial
applications is option pricing in these type of models (see
Nicolato and Venardos \cite{Nicolato} for some recent related results). These questions
will be addressed in future work.

\section*{Acknowledgements}
I would like to thank Professor D. Madan, Professor M. Fu, Professor E. Slud and Professor G. Skoulakis
for helpful discussions. Moreover, I would like to thank my
colleague Ziliang Li for his initial involvement in this project,
for bringing to my attention \cite{JMV}, for providing me with the
code of the Lambert-W function and for helpful suggestions.


\vspace{1cm}

\begin{verse}
\begin{small}
\hspace{3.2cm}{\sc Lefschetz Center for Dynamical Systems}\\
\hspace{3.2cm}{\sc Division of Applied Mathematics}\\
\hspace{3.2cm}{\sc Brown University} \\
\hspace{3.2cm}{\sc 182 George Street}\\
\hspace{3.2cm}{\sc Providence, RI, 02912, USA}\\
\hspace{3.2cm}{\sc E-mail: kspiliop@dam.brown.edu}
\end{small}
\end{verse}

\end{document}